\documentclass[12pt]{amsart}

\usepackage{fullpage}
\usepackage[utf8]{inputenc}
\usepackage[T1]{fontenc}
\usepackage{graphicx}
\usepackage{amsmath,amsfonts,amssymb,mathtools}
\usepackage{stmaryrd}
\usepackage{graphicx}
\usepackage{mathrsfs,dsfont}
\usepackage{amsmath,amsthm,amsthm}

\usepackage{hyperref}

\usepackage{enumerate}

\usepackage{color}

\newcommand{\E}{\mathbb{E}}
\newcommand{\R}{\mathbb{R}}
\newcommand{\N}{\mathbb{N}}

\newtheorem{theo}{Theorem}

\newtheorem{rem}{Remark}

\newtheorem{propo}{Proposition}

\newtheorem{defi}{Definition}

\newtheorem{hyp}{Assumption}

\begin{document}

\title{Influence of the regularity of the test functions for weak convergence in numerical discretization of  SPDE{\tiny s}}

\author{Charles-Edouard Br\'ehier}
\address{Univ Lyon, Université Claude Bernard Lyon 1, CNRS UMR 5208, Institut Camille Jordan, 43 blvd. du 11 novembre 1918, F-69622 Villeurbanne cedex, France}
\email{brehier@math.univ-lyon1.fr}

\begin{abstract}
This article investigates the role of the regularity of the test function when considering the weak error for standard discretizations of SPDEs of the form $dX(t)=AX(t)dt+F(X(t))dt+dW(t)$, driven by space-time white noise. In previous results, test functions are assumed (at least) of class $\mathcal{C}^2$ with bounded derivatives, and the weak order is twice the strong order.

We prove, in the case $F=0$, that to quantify the speed of convergence, it is crucial to control some derivatives of the test functions, even when the noise is non-degenerate. First, the supremum of the weak error over all bounded continuous functions, which are bounded by $1$, does not converge to $0$ as the discretization parameter vanishes. Second, when considering bounded Lipschitz test functions, the weak order of convergence is divided by $2$, {\it i.e.} it is not better than the strong order.

This is in contrast with the finite dimensional case, where the Euler-Maruyama discretization of elliptic SDEs $dY(t)=f(Y(t))dt+dB_t$ has weak order of convergence $1$ even for bounded continuous functions.
\end{abstract}


\date{}

\maketitle

\section{Introduction}

The numerical analysis of Stochastic Partial Differential Equations (SPDEs) has received a lot of attention in the last two decades, see for instance the recent monographs~\cite{JentzenKloeden:11},~\cite{Kruse:14} and~\cite{Lord_Powell_Shardlow:14}. Many temporal and spatial discretization schemes have been studied in the literature: Euler schemes, exponential Euler schemes, and spectral Galerkin methods, Finite Element methods.

In this article, we consider semilinear, parabolic, equations, with additive noise, of the type
$$
\left\{
\begin{array}{l}
dX= (\partial_{\xi\xi} X +F(X))dt+dW,\; t>0, \; \xi\in (0,1),\\
X(0,t)=X(1,t)=0,\\
X(\xi,0)=x(\xi),
\end{array}
\right.
$$
on the interval $(0,1)$, with homogeneous Dirichlet boundary conditions. More precisely, we consider Hilbert-space valued stochastic processes, which are solutions in $H=L^2(0,1)$ of
\begin{equation}\label{eq:SPDE_intro}
dX(t)=AX(t)dt+F(X(t))dt+dW(t),~X(0)=x,
\end{equation}
in the framework of~\cite{DPZ}, see equation~\eqref{eq:SPDE} and Section~\ref{sec:Setting} below for precise assumptions. The drift coefficient $F$ is assumed at least Lipschitz continuous, to ensure global well-posedness of mild solutions. In fact, we will mainly focus on the case $F=0$. The noise is given by a cylindrical Wiener process, which is a mathematical model for Gaussian space-time white noise.

We are interested in weak convergence rates for numerical approximations of $X(T)$, for arbitrary time $T\in(0,\infty)$. Recall that this notion corresponds to studying the weak error
\begin{equation}\label{eq:weak_error}
\E[\phi(X(T))]-\E[\phi(X_h(T))]
\end{equation}
where $X_h(T)$ is the numerical approximation of $X(T)$, obtained by temporal and/or spatial discretization of the equation (with discretization parameter $h\to 0$), and $\phi:H\to \R$ is a bounded continuous function. Recall also that strong convergence refers to the analysis of the strong error
\[
\E|X(T)-X_h(T)|.
\]
These notions have been extensively studied in the case of Stochastic Differential Equations (SDEs) of the type
\begin{equation}\label{eq:SDE_intro}
dY_t=f(Y_t)dt+\sigma(Y_t)dB_t,~Y_0=y\in\R^d,
\end{equation}
with smooth coefficients $f$ and $\sigma$, and a $d$-dimensional Brownian Motion $B$, see for instance the classic monographs~\cite{Kloeden_Platen:92},~\cite{Milstein_Tretyakov:04}.

Strong convergence for discretizations of the SPDEs~\eqref{eq:SPDE_intro}, also with multiplicative noise perturbation, have been studied, for instance, in~\cite{Davie_Gaines:01},~\cite{Gyongy_Millet:05},~\cite{JentzenKloeden:11},
~\cite{Printems:01},~\cite{Walsh:05},~\cite{Yan:05} (the list is not exhaustive). Results concerning weak convergence rates have essentially been obtained in the last decade, using different approaches. In the case of the stochastic equation with additive noise ($F=0$ in~\eqref{eq:SPDE_intro}), see~\cite{Debussche_Printems:09},
~\cite{Geissert_Kovacs_Larsson:09},
~\cite{Kovacs_Larsson_Lindgren:12},
~\cite{Kovacs_Larsson_Lindgren:13}. For semilinear equations, see
~\cite{Andersson_Larsson:16},
~\cite{Brehier_Debussche:17},
~\cite{Debussche:11},
~\cite{Tambue_Ngnotchouye:16},
~\cite{Wang:16},~\cite{Wang_Gan:13}, for an approach related to the Kolmogorov equation. See
~\cite{ConusJentzenKurniawan:14},
~\cite{Hefter_Jentzen_Kurniawan:16},
~\cite{JentzenKurniawan:15},
where a mild It\^o formula is used. Finally, for semilinear equations with additive noise, see~\cite{Andersson_Kruse_Larsson:16} and ~\cite{Brehier_Hairer_Stuart} for different approaches. Deriving weak convergence rates is fundamental in infinite dimension, see for instance~\cite{Lang:16}. Moreover, it is the appropriate notion for the approximation of invariant distribution (in the asymptotic regime $T\to \infty$), see~\cite{Brehier:14},~\cite{Brehier_Kopec:17},~\cite{Brehier_Vilmart:16}. The extension of the results of this article in this regime is straightforward.

\bigskip

The results in the references mentioned above can be roughly summarized as follows: if the strong error converges with order $r$, then the weak error converges with order $2r$, for functions $\phi$ which are sufficiently smooth, {\it i.e.} of class $\mathcal{C}^p$, bounded and with bounded derivatives of order $1,\ldots,p$, with $p\ge 2$ ($p$ depends on the model, for instance whether noise is additive or multiplicative):
\begin{equation}\label{eq:orders_intro}
\E|X(T)-X_h(T)|\le C(T)h^r~,\quad  \E[\phi(X(T))]-\E[\phi(X_h(T))]\le C(T)\|\phi\|_p h^{2r},
\end{equation}
where $\|\phi\|_p=\sup_{x\in H}|\phi(x)|+\sum_{j=1}^{p}\sup_{x\in H}|D^j\phi(x)|$. For spectral Galerkin discretization of the SPDE~\eqref{eq:SPDE_intro}, in dimension $N$, with $h=\frac{1}{N}$, one may choose $r\in[0,\frac{1}{2})$. For linear implicit Euler discretization of~\eqref{eq:SPDE_intro}, with time step size $h=\Delta t$, one may choose $r\in[0,\frac{1}{4})$.

This article investigates whether~\eqref{eq:orders_intro} holds true if $\|\phi\|_p$, where $p\ge 2$, is replaced with $\|\phi\|_1$ or $\|\phi\|_0$. This question is motivated by the positive answer for hypoelliptic SDEs, for instance in the additive noise case with constant $\sigma$; on the contrary, the contribution of this article shows that the answer is negative for SPDEs, and we exhibit some family of functions which allow us to identify the rate of convergence.

\bigskip

In the SDE case, consider the Euler-Maruyama discretization of~\eqref{eq:SDE_intro} (see~\eqref{eq:SDE_Euler}), with time step size $h$. Under an appropriate hypoellipticity assumption (which is satisfied in the additive non-degenerate noise case $\sigma(x)={\rm Id}$), using Malliavin calculus techniques and regularization effect in the associated Kolmogorov equation, the authors in~\cite{Bally_Talay_1:96},~\cite{Bally_Talay_2:96} (see also~\cite{Clement_Kohatsu-Higa_Lamberton:06}), have proved that the standard approach of~\cite{Talay:86}, to prove the weak error estimate for sufficiently regular functions,
\[
|\E[\phi(Y(T))]-\E[\phi(Y_h(T))]|\le C(T)\|\phi\|_p h
\]
with $p\ge 2$, can be extended with $\|\phi\|_0$ instead of $\|\phi\|_p$ on the right-hand side. In other words, weak convergence is also of order $1$ when considering bounded measurable test functions, in particular for bounded continuous test functions.

\bigskip

Our contribution is to prove that the situation is quite different for SPDEs. Note that thanks to~\eqref{eq:orders_intro}, the weak error~\eqref{eq:weak_error} converges to $0$ when $h\to 0$, for any given bounded continuous function $\phi$. Moreover, the Kolmogorov equation regularization effect also holds true in the infinite dimensional setting, thanks to non-degeneracy of the noise in~\eqref{eq:SPDE_intro}.

The main result of this paper states that the supremum over all bounded continuous functions, bounded by $1$, of the weak error~\eqref{eq:weak_error}, does not converge to $0$: the precise statement is Theorem~\ref{theo:0}. In addition, if one considers bounded Lipschitz continuous functions $\phi$, and set $\|\phi\|_1=\|\phi\|_0+\underset{x_1,x_2\in H}\sup\frac{|\phi(x_2)-\phi(x_1)|}{|x_2-x_1|}$, the weak error estimate in~\eqref{eq:orders_intro} is modified as
\[
\E[\phi(X(T))]-\E[\phi(X_h(T))]\le C(T)\|\phi\|_1 h^{r},
\]
see Theorem~\ref{theo:1} for a precise statement: there is a loss in the order of convergence. Equivalently, the optimal weak order for bounded Lipschitz continuous function is equal to the strong order for SPDEs, in the setting considered in this article.

The regularity of the test functions, and the control of derivatives, is thus essential to quantify the speed of convergence of the weak error~\eqref{eq:weak_error} for numerical discretization of SPDEs~\eqref{eq:SPDE_intro}.

\bigskip

Our proofs rely on academic examples of functions, which have low significance for concrete numerical approximation. It may be possible to define smaller families of non-regular test functions, for which uniform convergence of the numerical schemes holds true with better rates of convergence. This is expected to be obtained by the generalization of the finite dimensional approach of~\cite{Bally_Talay_1:96},~\cite{Bally_Talay_2:96}: regularization effect in the Kolmogorov equation and Malliavin calculus techniques. The identification of the appropriate setting is left for future works.

\bigskip

Why the regularity of the test functions matters for SPDEs may be explained by the properties of the solutions of associated Kolmogorov equations. Indeed, as emphasized in~\cite{Andersson_Larsson:16},
~\cite{Brehier_Debussche:17},
~\cite{Debussche:11},
 Sobolev-type regularity properties for the spatial derivatives of the solution of this infinite dimensional PDE are required to treat the most irregular terms in the error expansion. Similar arguments appear in~\cite{ConusJentzenKurniawan:14},
~\cite{Debussche_Printems:09} and related articles. The regularity estimates have singularities at the initial time, even when the test function (seen as the initial condition of the Kolmogorov equation) is regular.

For SDEs, the Kolmogorov equation preserves regularity of the initial condition. Singularities only appear when a regularization effect is needed, in an hypoelliptic setting.

For SPDEs, exhibiting a rate of convergence in the error analysis is only possible when using some spatial regularity property, as mentioned above. The better the spatial regularity, the greater the order of convergence, but the stronger the singularity -- with the constraint of remaining integrable. This approach yields the optimal order of convergence for regular test functions. Weakening the regularity condition on the test functions then introduces even stronger singularities, and less spatial regularity may be used: in turn the order of convergence decreases. The optimality of these heuristic arguments is validated by Theorems~\ref{theo:0} and~\ref{theo:1}.

\bigskip

The article is organized as follows. Assumptions on the model and numerical discretization schemes are introduced in Section~\ref{sec:Setting}. Section~\ref{sec:auxiliary} describes important spatial regularity properties, which are very different for the discretized versions, compared with the exact solution. Our main results, Theorem~\ref{theo:0} (bounded continuous functions) and Theorem~\ref{theo:1} are stated in Section~\ref{sec:results}. Detailed proofs are provided in Section~\ref{sec:proofs}.

\section{Setting}\label{sec:Setting}

\subsection{Model and assumptions}

The model in this article is given by a Stochastic Partial Differential Equation (SPDE),
\begin{equation}\label{eq:SPDE}
dX(t)=AX(t)dt+dW(t),\quad X(0)=0,
\end{equation}
{\it i.e.} by Equation~\ref{eq:SPDE_intro}, with $F=0$. This choice is sufficient for our purpose and does not change the conclusions of this article.


The initial condition in~\eqref{eq:SPDE} is set to $0$ for simplicity. Extending the results of this article to arbitrary initial conditions is straightforward.

\subsubsection{Linear operator $A$}

Denote by $\langle\cdot,\cdot\rangle$, resp. $|\cdot|$, the inner product, resp. the norm, in the separable Hilbert space $H=L^2(0,1)$.

The operator $A$ in the SPDE satisfies the following conditions.
\begin{hyp}\label{ass:A}
The mapping $A$ is an unbounded, self-adjoint, linear operator on $H$.

Define, for all $n\in\N=\left\{1,\ldots\right\}$,
\[
\lambda_n=\pi^2n^2~,\quad e_n=\sqrt{2}\sin\bigl(n\pi \cdot\bigr).
\]
Then the operator $A$ and its domain $D(A)$ are given by
\[
Ax=\sum_{n\in\N}-\lambda_n \langle x,e_n\rangle~,\quad\forall~x\in D(A)=\left\{x\in H~;~\sum_{n\in\N}\lambda_n^2\langle x,e_n\rangle^2<\infty\right\}.
\]
\end{hyp}

Recall that $\bigl(e_n\bigr)_{n\in\N}$ is a complete orthonormal system of $H$.

Introduce the following notation.
\begin{defi}
\begin{enumerate}
\item The operator $A$ generates a strongly-continuous semigroup $\bigl(e^{tA}\bigr)_{t\ge 0}$ on $H$, with
\[
e^{tA}x=\sum_{n\in\N}e^{-\lambda_n t}\langle x,e_n\rangle~,\quad\forall~x\in H,t\ge 0.
\]
\item For all $\alpha\in[0,1]$, set
\[
|x|_\alpha=\bigl(\sum_{n\in\N}\lambda_n^{2\alpha}\langle x,e_n\rangle^2\bigr)^{\frac{1}{2}}\in[0,\infty],~\quad \forall~x\in H.
\]
\end{enumerate}
\end{defi}

\subsubsection{Cylindrical Wiener process}

\begin{hyp}
Let $\bigl(\Omega,\mathcal{F},\mathbb{P}\bigr)$ denote a probability space, expectation is denoted by $\E$.

Let $\bigl(\beta_n\bigr)_{n\in\N}$ be a sequence of independent standard $\mathbb{R}$-valued Wiener processes.

Then set, for all $t\ge 0$,
\begin{equation}\label{eq:W}
W(t)=\sum_{n\in\N}\beta_n(t)e_n.
\end{equation}
\end{hyp}

It is a standard fact that, for all $t\ge 0$, almost surely the series in~\eqref{eq:W} does not converge in $H$. However, if $\Phi\in\mathcal{L}(H)$ is an Hilbert-Schmidt operator, then $\Phi W(t)=\sum_{n\in\N}\beta_n(t) \Phi e_n$ is a Wiener process in $H$, with covariance operator $\Phi\Phi^{\star}$.

\subsubsection{Mild solution}

Solutions of the SPDE~\eqref{eq:SPDE} are interpreted in the mild sense: the unique solution, which is often called stochastic convolution, is given by
\[
X(t)=\int_{0}^{t}e^{(t-s)A}dW(s)=\sum_{n\in\N}\bigl(\int_{0}^{t}e^{-\lambda_n(t-s)}d\beta_n(s)\bigr)e_n, \quad t\in(0,\infty).
\]

With the notation $X_n(t)=\langle X(t),e_n\rangle=\int_{0}^{t}e^{-\lambda_n(t-s)}d\beta_n(s)$, the process $\bigl(X_n(t)\bigr)_{t\ge 0}$ are independent Ornstein-Uhlenbeck processes. Thus, $\bigl(X(t)\bigr)_{t\ge 0}$ is a centered Gaussian process with values in $H$. Let $\mu_t$ denote the law of $X(t)$, {\it i.e.} the centered Gaussian probability distribution on $H$ with covariance operator $Q_t\in\mathcal{L}(H)$, given by $Q_t e_n=\frac{1}{2\lambda_n}\bigl(1-e^{-2\lambda_n t}\bigr)=\E\bigl[|X_n(t)|^2\bigr]$ for all $n\in\N$ and $t\ge 0$. 

\subsection{Numerical schemes}

Space and time discretization schemes are defined below. One may also consider full-discretization schemes obtained by combining these two procedures.

\subsubsection{Space discretization: spectral Galerkin method}

For every $N\in\N$, let $P_N\in\mathcal{L}(H)$ denote the orthogonal projection onto the finite-dimensional subspace ${\rm Span}\bigl(e_1,\ldots,e_N\bigr)$:
\[
P_Nx=\sum_{n=1}^{N}\langle x,e_n\rangle e_n,\quad \forall~x\in H.
\]

The process $X^{(N)}$ obtained by discretization in space of the SPDE~\eqref{eq:SPDE}, is solution of
\[
dX^{(N)}(t)=AX^{(N)}(t)dt+P_NdW(t)~,\quad X^{(N)}(0)=0.
\]
In fact, $X^{(N)}(t)=P_NX(t)$, for all $t\ge 0$ and $N\in\N$.

Let then $\mu_t^{(N)}$ denote the law of the random variable $X^{(N)}$: it is a centered Gaussian probability distribution, with covariance operator $P_NQ_t(P_N)^{\star}=P_NQ_t$.

\subsubsection{Time discretization: linear implicit Euler scheme}

Let $\Delta t>0$ denote a time-step size, without restriction we assume $\Delta t\in(0,1)$. The scheme is defined such that for all $k\in\N_0=\left\{0,1,\ldots\right\}$,
\[
X_{k+1}^{\Delta t}=X_k^{\Delta t}+\Delta t AX_{k+1}^{\Delta t}+\Delta W_k~,\quad X_0^{\Delta t}=0,
\]
with Wiener increments $\Delta W_k=W\bigl((k+1)\Delta t\bigr)-W\bigl(k\Delta t\bigr)$.

Rigorously,
\[
X_{k+1}^{\Delta t}=S_{\Delta t}X_k+S_{\Delta t}\Delta W_k,
\]
where $S_{\Delta t}=\bigl(I-\Delta t A\bigr)^{-1}$ is a linear, self-adjoint, Hilbert-Schmidt, operator on $H$.

As a consequence, for every $k\in\N$,
\[
X_k^{\Delta t}=\sum_{\ell=0}^{k-1}S_{\Delta t}^{k-\ell}\Delta W_\ell,
\]
and the law $\nu_k^{\Delta t}$ of $X_k^{\Delta t}$ is a centered Gaussian probability distribution, with covariance operator
\[
Q_k^{\Delta t}=\Delta t\sum_{\ell=0}^{k-1}S_{\Delta t}^{2(k-\ell)}.
\]

\subsection{Space regularity properties}\label{sec:auxiliary}

The aim of this section is to provide some important results concerning the moments $\int_{H}|x|_\alpha^2 \mu(dx)$, for different values of $\alpha\in[0,1]$. The parameter $\alpha$ is interpreted as indicating space regularity of the process. We emphasize on the key observation: the behaviors are different when considering, on the one hand, $\mu=\mu_t$, and, on the other hand, $\mu=\mu_t^{(N)}$ or $\mu=\nu_n^{\Delta t}$, which are obtained by the discretization schemes. We will take advantage of this property in the study of the orders of convergence for bounded continuous test functions.

First, consider the law $\mu_t$ at time $t$, of the solution of the SPDE~\eqref{eq:SPDE}: for any $t\in(0,\infty)$,
\begin{equation}\label{eq:norm_alpha_exact}
\int_{H}|x|_\alpha^2\mu_t(dx)=\sum_{n\in\N}\frac{1}{2\lambda_n^{1-2\alpha}}\bigl(1-e^{-2\lambda_n t}\bigr)<\infty \quad \Longleftrightarrow\quad \alpha\in[0,\frac{1}{4}).
\end{equation}

Now, consider the law $\mu_{t}^{(N)}$, at time $t$, obtained by spatial discretization: for every $t\in[0,\infty)$,
\begin{equation}\label{eq:norm_alpha_space}
\begin{gathered}
\int_{H}|x|_{\alpha}^{2}\mu_t^{(N)}(dx)<\infty,\quad \forall~\alpha\in[0,1],~\forall~N\in\N,\\
\underset{N\in\N}\sup \int_{H}|x|_{\alpha}^{2}\mu_t^{(N)}(dx)<\infty \quad \Longleftrightarrow\quad \alpha\in[0,\frac{1}{4}).
\end{gathered}
\end{equation}

Finally, consider the law $\nu_{k}^{\Delta t}$, at time $k$, obtained by the temporal discretization: for every $k\in\N$
\begin{equation}\label{eq:norm_alpha_time}
\begin{gathered}
\int_H |x|_\alpha^2 \nu_k^{\Delta t}(dx)<\infty,\quad \forall~\alpha\in[0,\frac{3}{4}),~\forall~\Delta t\in(0,1),\\
\underset{\Delta t\in (0,1)}\sup \int_H |x|_\alpha^2 \nu_k^{\Delta t}(dx)<\infty \quad \Longleftrightarrow \quad \alpha\in[0,\frac{1}{4}).
\end{gathered}
\end{equation}

Observe that in~\eqref{eq:norm_alpha_space} and~\eqref{eq:norm_alpha_time}, one recovers the same behavior as in~\eqref{eq:norm_alpha_exact}, only when the supremum over all discretization parameters ($N\in\N$ and $\Delta t\in(0,1)$) is computed. For fixed values of these parameters, some larger values of $\alpha\ge \frac{1}{4}$ are allowed.

The proofs of estimates in~\eqref{eq:norm_alpha_exact} and~\eqref{eq:norm_alpha_space} are straightforward. For completeness, let us give a detailed proof of the estimates in~\eqref{eq:norm_alpha_time}. Similar arguments will be used again below.

To prove the first statement in~\eqref{eq:norm_alpha_time}, let $\Delta t>0$, $k\in\N$, and $\alpha\in[0,1]$, then
\begin{align*}
\int_H |x|_\alpha^2 \nu_k^{\Delta t}(dx)&=\sum_{n\in\N}\lambda_n^{2\alpha}\langle Q_k^{\Delta t}e_n,e_n\rangle=\Delta t\sum_{n\in\N}\lambda_n^{2\alpha}\sum_{\ell=1}^{k}\frac{1}{(1+\lambda_n\Delta t)^{2\ell}}\\
&=\sum_{n\in\N}\frac{\lambda_n^{2\alpha}}{\lambda_n(2+\lambda_n\Delta t)}\bigl(1-\frac{1}{(1+\lambda_n\Delta t)^{2k}}\bigr)\\
&<\infty \quad \Longleftrightarrow \quad \alpha\in[0,\frac{3}{4}).
\end{align*}
To prove the second statement, first assume $\alpha\in[0,\frac{1}{4})$, then $1-2\alpha>\frac{1}{2}$, thus for all $\Delta t\in(0,1)$, and all $k\in\N$,
\[
\int_H |x|_\alpha^2 \nu_k^{\Delta t}(dx)\le \sum_{n\in\N}\frac{1}{2\lambda_{n}^{1-2\alpha}}<\infty.
\]
Now assume that $\alpha\ge \frac{1}{4}$. By a monotonicity argument, it is sufficient to consider the case $\alpha=\frac{1}{4}$. Let $M\in\N$ be an auxiliary integer, and choose $\Delta t=\frac{1}{N^2}$, with $N\in\N$, $N\ge M$.
\begin{align*}
\int_H |x|_{\frac{1}{4}}^2 \nu_k^{\frac{1}{N^2}}(dx)&= \sum_{n\in\N}\frac{1}{\pi n(2+\pi^2\frac{n^2}{N^2})}\bigl(1-\frac{1}{(1+\pi^2\frac{n^2}{N^2})^{2k}}\bigr)\\
&\ge \frac{1}{N}\sum_{n\ge \frac{N}{M}}\frac{1}{\pi \frac{n}{N}(2+\pi^2\frac{n^2}{N^2})}\bigl(1-\frac{1}{(1+\pi^2\frac{n^2}{N^2})^{2k}}\bigr)\\
&\underset{N\to \infty}\to \frac{1}{\pi}\int_{\frac{\pi}{M}}^{\infty}\frac{1}{z(2+z^2)}(1-\frac{1}{(1+z^2)^{2k}})dz,
\end{align*}
by a Riemann sum argument. Then $\underset{N\to\infty}\liminf \int_H |x|_{\frac{1}{4}}^2 \nu_k^{\frac{1}{N^2}}(dx)\ge \frac{1}{\pi}\int_{0}^{\infty}\frac{1}{z(2+z^2)}(1-\frac{1}{(1+z^2)^{2k}})dz=\infty$, taking $M\to\infty$. This concludes the proof of the equivalence statement in~\eqref{eq:norm_alpha_time}.

\section{Main results}\label{sec:results}

Introduce the following notation:
\begin{itemize}
\item $\|\phi\|_0=\underset{x\in H}\sup|\phi(x)|$, for $\phi\in\mathcal{C}^{0}(H,\R)$, bounded and continuous functions from $H$ to $\R$,
\item $\|\phi\|_1=\|\phi\|_0+\underset{x,y\in H, x\neq y}\sup\frac{|\phi(y)-\phi(x)|}{|y-x|}$, for $\phi\in\mathcal{C}^{0,1}(H,\R)$, bounded and Lipschitz continuous functions from $H$ to $\R$,
\item $\|\phi\|_2=\underset{x\in H}\sup|\phi(x)|+\underset{x\in H,h\in H,|h|\le 1}\sup|D\phi(x).h|+\underset{x\in H,h_1,h_2\in H,|h_1|\le 1,|h_2|\le 1}\sup|D^2\phi(x).(h_1,h_2)|$, for $\phi\in\mathcal{C}^{2}(H,\R)$, bounded functions from $H$ to $\R$ of class $\mathcal{C}^2$, with bounded first and second order derivatives.
\end{itemize}

\subsection{Statements}

The main result of this article is Theorem~\ref{theo:0}, which may be interpreted as follows: there is no rate of convergence to $0$, for the weak error, when considering the supremum over all bounded and continuous functions.
\begin{theo}\label{theo:0}
Let $T\in(0,\infty)$. Then
\begin{equation}\label{eq:theo_0}
\begin{gathered}
\underset{N\to \infty}\limsup \underset{\phi\in\mathcal{C}^0(H,\R),\|\phi\|_{0}\le 1}\sup |\int \phi d\mu_T-\int \phi d\mu_{T}^{(N)}|>0,\\
\underset{\Delta t\to 0}\limsup \underset{\phi\in\mathcal{C}^0(H,\R),\|\phi\|_{0}\le 1}\sup |\int \phi d\mu_T-\int \phi d\nu_{\lfloor\frac{T}{\Delta t}\rfloor}^{\Delta t}|>0.
\end{gathered}
\end{equation}
\end{theo}

The proof of Theorem~\ref{theo:0} is postponed to Section~\ref{sec:proof_0}.

To explain why the statement of Theorem~\ref{theo:0} may be surprising, recall that strong convergence results, with order in $[0,\frac{1}{4})$, are available: for every $r\in[0,\frac{1}{4})$, and every $T\in(0,\infty)$,
\begin{equation}\label{eq:strong}
\limsup_{N\to \infty}\lambda_{N}^{2r}\E\big|X(T)-X^{(N)}(T)\big|^2<\infty~,\quad
\limsup_{\Delta t\to 0}\frac{1}{\Delta t^{2r}}\E\big|X(T)-X_{\lfloor \frac{T}{\Delta t}\rfloor}^{\Delta t}\big|^2<\infty.
\end{equation}
Thus, for any bounded and continuous function $\phi\in\mathcal{C}^0(H,\R)$, the convergence below is valid:
\[
\int \phi d\mu_{T}^{(N)}\underset{N\to \infty}\to \int \phi d\mu_T \quad,\quad \int \phi d\nu_{\lfloor\frac{T}{\Delta t}\rfloor}^{\Delta t}\underset{\Delta t\to 0}\to\int \phi d\mu_T.
\]
However, the supremum of the error over all bounded continuous functions, bounded by $1$, does not converge to $0$.

As will become clear in the proof of Theorem~\ref{theo:0}, see the stronger statement~\eqref{eq:theo_0_bis} below, the issue is not the regularity of the functions $\phi$ -- smooth functions are used -- but the lack of control of the growth of the derivatives.

It is also worth mentioning that if one considers the set of bounded measurable test functions, instead of continuous test functions, in~\eqref{eq:theo_0}, the result is straightforward, see Remark~\ref{rem:rem}. Indeed, this corresponds to looking at the total variation distance between $\mu_T$, and $\mu_T^{(N)}$ of $\nu_{\lfloor\frac{T}{\Delta t}\rfloor}^{\Delta t}$, and due to the results of Section~\ref{sec:auxiliary}, these distributions are singular.

We also prove the following statement, Theorem~\ref{theo:1}, which may be interpreted as follows: the best order of convergence, for the weak error, when considering the supremum over all bounded and Lipschitz continuous functions, is equal to the strong order of convergence. 
\begin{theo}\label{theo:1}
Let $T\in(0,\infty)$. Then
\begin{equation}\label{eq:theo_1}
\begin{gathered}
\underset{N\to \infty}\limsup~ \lambda_{N}^{r}\underset{\phi\in\mathcal{C}^{0,1}(H,\R),\|\phi\|_{1}\le 1}\sup |\int \phi d\mu_T-\int \phi d\mu_{T}^{(N)}|=\begin{cases} 0,~\forall~r\in[0,\frac{1}{4})\\
\infty,~\forall~r\in(\frac{1}{4},\frac{1}{2})\end{cases},
\\
\underset{\Delta t\to 0}\limsup~\frac {1}{\Delta t^r}\underset{\phi\in\mathcal{C}^{0,1}(H,\R),\|\phi\|_{1}\le 1}\sup |\int \phi d\mu_T-\int \phi dd\nu_{\lfloor\frac{T}{\Delta t}\rfloor}^{\Delta t}|=\begin{cases} 0,~\forall~r\in[0,\frac{1}{4})\\
\infty,~\forall~r\in(\frac{1}{4},\frac{1}{2})\end{cases}.
\end{gathered}
\end{equation}
\end{theo}

The results in Theorem~\ref{theo:1} in the regime $r\in[0,\frac{1}{4})$ are not new, they are straightforward applications of the strong convergence estimates in~\eqref{eq:strong}. The case $r\in(\frac{1}{4},\frac{1}{2})$ is treated in Section~\ref{sec:proof_1}.

For comparison, we state an additional result, considering test functions of class $\mathcal{C}^2$, bounded, and with bounded first and second order derivatives.
\begin{propo}\label{propo:2}
Let $T\in(0,\infty)$.
\begin{equation}\label{eq:propo_2}
\begin{gathered}
\underset{N\to \infty}\limsup~ \lambda_{N}^{r}\underset{\phi\in\mathcal{C}^{2}(H,\R),\|\phi\|_{2}\le 1}\sup |\int \phi d\mu_T-\int \phi d\mu_{T}^{(N)}|=\begin{cases} 0,~\forall~r\in[0,\frac{1}{2})\\
\infty,~\forall~r\in(\frac{1}{2},1)\end{cases},
\\
\underset{\Delta t\to 0}\limsup~\frac {1}{\Delta t^r}\underset{\phi\in\mathcal{C}^{2}(H,\R),\|\phi\|_{2}\le 1}\sup |\int \phi d\mu_T-\int \phi dd\nu_{\lfloor\frac{T}{\Delta t}\rfloor}^{\Delta t}|=\begin{cases} 0,~\forall~r\in[0,\frac{1}{2})\\
\infty,~\forall~r\in(\frac{1}{2},1)\end{cases}.
\end{gathered}
\end{equation}
\end{propo}

The result of Proposition~\ref{propo:2}, in the regime $r\in[0,\frac{1}{2})$, has been proved in a more general setting, for semilinear versions of~\eqref{eq:SPDE}, see the references in the introduction. In the case of multiplicative noise, the results require that $\phi$ is at least of class $\mathcal{C}^3$, however the order of convergence remains equal to $\frac{1}{2}$ for such test functions. The case $r\in(\frac{1}{2},1)$ is obtained using the lower bounds from~\cite{ConusJentzenKurniawan:14}.

Note that Theorems~\ref{theo:0} and~\ref{theo:1} are also valid when looking at the regime $T\to\infty$, {\it i.e.} at the level of the invariant distributions of the process and of its discretized versions.

Comparing Theorems~\ref{theo:0},~\ref{theo:1} and Proposition~\ref{propo:2} reveals that in infinite dimension, regularity of the test functions and control of derivatives plays an important role in the analysis of the numerical error in the weak sense.

\subsection{Comparison with the finite dimensional situation}

The situation described by Theorems~\ref{theo:0} and~\ref{theo:1}, and Proposition~\ref{propo:2}, is specific to the infinite dimensional situation. Indeed, when considering Euler-Maruyama discretization of hypoelliptic SDEs (in finite dimension), the order of convergence (equal to $1$) does not change when considering either bounded continuous functions, or bounded and Lipschitz continuous functions, or functions of class $\mathcal{C}^2$.

Indeed, consider a SDE in $\mathbb{R}^d$ (see Equation~\eqref{eq:SDE_intro}, with additive non-degenerate noise),
\begin{equation}\label{eq:SDE}
dY(t)=f(Y(t))dt+dB_t, Y(0)=y_0,
\end{equation}
where $\bigl(B_t\bigr)_{t\ge 0}$ is a $d$-dimensional standard Wiener process, and $f:\R^d\to \R^d$ is a smooth bounded function, with bounded derivatives.

Consider its Euler-Maruyma discretization, with time step size $\Delta t>0$: for $k\in\N_0$,
\begin{equation}\label{eq:SDE_Euler}
Y_{k+1}^{\Delta t}=Y_k^{\Delta t}+\Delta tf(Y_k^{\Delta t})+B((k+1)\Delta t)-B(k\Delta t)~,\quad Y_0^{\Delta t}=y_0.
\end{equation}
The strong order of convergence in this case is equal to $1$ (this is due to the fact that the noise is additive, it would be equal to $\frac{1}{2}$ in general):
\[
\underset{\Delta t\to 0}\limsup \frac{1}{\Delta t^2}\E|Y(T)-Y_{\lfloor \frac{T}{\Delta t}\rfloor}^{\Delta t}|^2\in(0,\infty).
\]

Then, it is a remarkable fact that when considering bounded measurable test functions, one still obtains an error which is of order $1$, see~\cite{Bally_Talay_1:96},~\cite{Bally_Talay_2:96},
\[
\underset{\Delta t\to 0}\limsup\frac{1}{\Delta t}\underset{\phi\in\mathcal{C}^0(H,\R),\|\phi\|_{0}\le 1}\sup |\E\phi(Y(T))-\E\phi(Y_{\lfloor \frac{T}{\Delta t}\rfloor}^{\Delta t})| \in(0,\infty),
\]
for every $T\in(0,\infty)$. Equivalently, there exists $C(T)\in(0,\infty)$, such that for every bounded continuous function $\phi$,
\begin{equation}\label{eq:C0_SDE}
|\E\phi(Y(T))-\E\phi(Y_{\lfloor \frac{T}{\Delta t}\rfloor}^{\Delta t})|\le C(T)\|\phi\|_0\Delta t.
\end{equation}

Theorem~\ref{theo:0} indicates that in infinite dimension, the generalization of~\eqref{eq:C0_SDE} is not valid, both for the standard and widely used time and space discretization schemes we have considered.

\section{Proofs}\label{sec:proofs}

\subsection{Bounded continuous test functions: proof of Theorem~\ref{theo:0}}\label{sec:proof_0}

In fact, a slightly stronger result than Theorem~\ref{theo:0} is proved below:
\begin{equation}\label{eq:theo_0_bis}
\underset{N\to \infty}\limsup~\underset{\phi\in\Phi}\sup |\int \phi d\mu_T-\int \phi d\mu_{T}^{(N)}|>0~,\quad
\underset{\Delta t\to 0}\limsup~\underset{\phi\in\Phi}\sup |\int \phi d\mu_T-\int \phi d\nu_{\lfloor\frac{T}{\Delta t}\rfloor}^{\Delta t}|>0,
\end{equation}
where $\Phi\subset\mathcal{C}^\infty(H,\R)$ is such that $\|\phi\|_0\le 1$ for all $\phi\in\Phi$. In the examples given below, the functions $\phi$ are smooth and have bounded derivatives of any order, however only $\|\phi\|_0$ is uniformly bounded over $\Phi$ -- precisely, $\sup\left\{\|\phi\|_1,~\phi\in\Phi \right\}=\infty$.

We provide two different examples of sets $\Phi$. The first family is constructed using the results of Section~\ref{sec:auxiliary}, concerning regularity properties of the discretized versions of the SPDE, see Remark~\ref{rem:rem}. The second family contains functions with arbitrarily fast oscillations, and is treated using some Riemann sums arguments. This proof is instructive, similar arguments appear for proving Theorem~\ref{theo:1}.


\subsubsection{First proof}

Define $\Phi^1=\left\{\phi_{\epsilon,M}^1,~\epsilon\in(0,1),M\in\N
\right\}$, where
\begin{equation}
\phi_{\epsilon,M}^1(x)=\exp\bigl(-\epsilon |P_Mx|_{\frac{1}{4}}^{2}\bigr),\quad\forall~x\in H.
\end{equation}
Then $\phi_{\epsilon,M}^2\in\mathcal{C}^\infty(H,\R)$, and $\|\phi_{\epsilon,M}^{1}\|_{0}=1$. However, $\sup\left\{\|\phi\|_1,~\phi\in\Phi^1\right\}=\infty$.

Let
\[
\begin{gathered}
\delta_{1}^{1}(N)=\underset{\phi\in\Phi^1}\sup |\int \phi d\mu_T-\int \phi d\mu_{T}^{(N)}|,\\
\delta_{2}^{1}(\Delta t)=\underset{\phi\in\Phi^1}\sup |\int \phi d\mu_T-\int \phi d\nu_{\lfloor\frac{T}{\Delta t}\rfloor}^{\Delta t}|.
\end{gathered}
\]

For every $N\in\N$, $\epsilon\in(0,1)$, letting $M\to\infty$ gives
\[
\delta_1^1(N)\ge \big|\E[e^{-\epsilon|P_MX(T)|_{\frac{1}{4}}^{2}}]-\E[e^{-\epsilon|P_MX^{(N)}(T)|_{\frac{1}{4}}^{2}}]\big|\ge \big|0-\E[e^{-\epsilon|X^{(N)}(T)|_{\frac{1}{4}}^{2}}]|,
\]
where almost surely $|X^{(N)}(T)|_{\frac{1}{4}}^{2}<\infty$, thanks to~\eqref{eq:norm_alpha_space}. On the contrary, using~\eqref{eq:norm_alpha_exact}, $\E|X(T)|_{\frac{1}{4}}^{2}=\infty$, and in fact almost surely $|X(T)|_{\frac{1}{4}}^{2}=\infty$. More precisely,
\begin{align*}
\E[e^{-\epsilon|P_MX(T)|_{\frac{1}{4}}^{2}}]&=\prod_{m=1}^{M}\E[e^{-\epsilon \lambda_m^{\frac{1}{2}}|\langle X(T),e_m\rangle|^2}]\\
&=\prod_{m=1}^{M}\bigl(1+\frac{\epsilon\sqrt{\lambda_m}}{\lambda_m}(1-e^{-2\lambda_mT})\bigr)^{-\frac{1}{2}}\\
&=\exp\bigl(-\frac{1}{2}\sum_{m=1}^{M}\log\bigl(1+\frac{\epsilon\sqrt{\lambda_m}}{\lambda_m}(1-e^{-2\lambda_mT})\bigr)\bigr)\\
&\underset{M\to \infty}\to 0.
\end{align*}

Similarly,
\[
\delta_{2}^{1}(\Delta t)\ge \E[e^{-\epsilon|X_{\lfloor\frac{T}{\Delta t}\rfloor}^{\Delta t}|_{\frac{1}{4}}^{2}}],
\]
with $|X_{\lfloor\frac{T}{\Delta t}\rfloor}^{\Delta t}|_{\frac{1}{4}}^{2}<\infty$ almost surely, thanks to~\eqref{eq:norm_alpha_time}.

Finally, letting $\epsilon\to 0$, for all $N\in\N$ and $\Delta t\in (0,1)$
\[
\delta_{1}^{1}(N)\ge 1~,\quad \delta_{2}^{1}(\Delta t)\ge 1.
\]
Thus
\begin{gather*}
\underset{N\to \infty}\limsup~\underset{\phi\in\Phi^1}\sup |\int \phi d\mu_T-\int \phi d\mu_{T}^{(N)}|=\underset{N\to \infty}\limsup~\delta_{1}^{1}(N)\ge 1,\\
\underset{\Delta t\to 0}\limsup~\underset{\phi\in\Phi^1}\sup |\int \phi d\mu_T-\int \phi d\nu_{\lfloor\frac{T}{\Delta t}\rfloor}^{\Delta t}|=\underset{\Delta t\to 0}\limsup~\delta_{2}^{1}(\Delta t)\ge 1,
\end{gather*}
hence~\eqref{eq:theo_0}. This concludes the first proof of Theorem~\ref{theo:0}.

\begin{rem}\label{rem:rem}
Consider the bounded measurable function $\overline{\phi^1}:H\to \R$, given by
\[
\overline{\phi^1}(x)=\mathds{1}_{|x|_{\frac{1}{4}}<\infty}.
\]
Then, thanks to the regularity results~\eqref{eq:norm_alpha_exact},~\eqref{eq:norm_alpha_space} and~\eqref{eq:norm_alpha_time} for all $N\in\N$ and $\Delta t\in(0,1)$,
\[
\int \overline{\phi^1} d\mu_T=0~,\quad \int \overline{\phi^1}  d\mu_{T}^{(N)}=\int\overline{\phi^1} d\nu_{\lfloor\frac{T}{\Delta t}\rfloor}^{\Delta t}=1.
\]
This means that, $\mu_{T}$ and $\mu_{T}^{(N)}$, resp. $\mu_{T}$ and $\nu_{\lfloor\frac{T}{\Delta t}\rfloor}^{\Delta t}$, are singular probability distributions on the infinite dimensional space $H$, and their distance in total variation is equal to $1$.

When $\epsilon\to 0$ and $M\to\infty$, the continuous functions $\phi_{\epsilon,M}^1$ converge pointwise to the measurable function $\overline{\phi}^1$, hence the idea of the proof.
\end{rem}

\subsubsection{Second proof}

Define $\Phi^2=\left\{\phi_M^2,~M\in\N\right\}$, where
\[
\phi_M^2(x)=\exp\bigl(i\sqrt{M}\langle \theta_M,x\rangle\bigr),~\theta_M=\sum_{m=\frac{M}{2}}^{M}e_m,\quad\forall~x\in H.
\]
In this example, it is convenient to consider complex-valued functions, however it is straightforward to get rid of this issue.

Like above, $\Phi^2\subset\mathcal{C}^\infty(H,\mathbb{C})$, $\|\phi\|_0=1$ for all $\phi\in\Phi^2$, and $\sup\left\{\|\phi\|_1,~\phi\in\Phi^2\right\}=\infty$. Let
\[
\delta_{1}^{2}(N)=\underset{\phi\in\Phi^2}\sup |\int \phi d\mu_T-\int \phi d\mu_{T}^{(N)}|~,\quad
\delta_{2}^{2}(\Delta t)=\underset{\phi\in\Phi^2}\sup |\int \phi d\mu_T-\int \phi d\nu_{\lfloor\frac{T}{\Delta t}\rfloor}^{\Delta t}|.
\]

First, focus on $\delta_{1}^{2}(N)$. Observe that for $M\ge 2N+1$, $\langle \theta_M,X^{(N)}(T)\rangle=0$ almost surely, hence $\int\phi_M^2 d\mu_T^{(N)}=1$. On the contrary,
\begin{align*}
\int\phi_M^2 d\mu_T&=\E[e^{i\sqrt{M}\langle \theta_M,X(T)\rangle}]=\exp\bigl(-M\sum_{m=\frac{M}{2}}^{M}\frac{1}{2\lambda_m}(1-e^{-2\lambda_m T})\bigr)\\
&=\exp\bigl(-\frac{1}{M}\sum_{m=\frac{M}{2}}^{M}\frac{1}{2\pi^2\bigl(\frac{m}{M}\bigr)^2}+{\rm o}(1)\bigr)\\
&\underset{M\to\infty}\to \exp\bigl(-\int_{\frac{1}{2}}^{1}\frac{1}{2\pi^2 z^2}dz\bigr),
\end{align*}
using a Riemann sum argument, and $M\sum_{m=\frac{M}{2}}^{M}\frac{1}{2\lambda_m}e^{-2\lambda_m T}={\rm O}\bigl(e^{-2\lambda_{\frac{M}{2}} T} \bigr)\underset{M\to \infty}\to 0$.


As a consequence, for all $N\in\N$, and letting $M\to\infty$ (with $M\ge 2N+1$), one obtains
\[
\delta_{1}^{2}(N)\ge 1-\exp\bigl(-\int_{\frac{1}{2}}^{1}\frac{1}{2\pi^2 z^2}dz\bigr)>0.
\]

Second, focus on $\delta_{2}^{2}(\Delta t)$. In order to use a Riemann sum argument, it is convenient to choose $\Delta t=\frac{T}{M^2}$, and to write
\[
\delta_{2}^{2}\bigl(\frac{T}{M^2}\bigr)\ge |\int \phi_M^2 d\mu_T-\int\phi_M^2d\nu_{M^2}^{\Delta t}|,
\]
with $\int\phi_M^2 d\mu_T\underset{M\to \infty}\to \exp\bigl(-\int_{\frac{1}{2}}^{1}\frac{1}{2\pi^2 z^2}dz\bigr)$, as above, and
\begin{align*}
\int\phi_M^2d\nu_{M^2}^{\Delta t}&=\E[\exp(i\sqrt{M}\langle \theta_M,X_{M^2}^{\Delta t}\rangle\bigr)\\
&=\exp\bigl(-M\sum_{m=\frac{M}{2}}^{M}\frac{1}{\lambda_m(2+\lambda_m\Delta t)}\bigl(1-\frac{1}{(1+\lambda_m\Delta t)^{2M^2}}\bigr)\\
&=\exp\bigl(-\frac{1}{M}\sum_{m=\frac{M}{2}}^{M}\frac{1}{\pi^2(\frac{m}{M})^2(2+\pi^2(\frac{m}{M})^2)}+{\rm o}(1)\bigr)\\
&\underset{M\to \infty}\to \exp\bigl(-\int_{\frac{1}{2}}^{1}\frac{1}{\pi^2z^2(2+\pi^2z^2)}dz\bigr),
\end{align*}
using a Riemann sum argument, and $M\sum_{m=\frac{M}{2}}^{M}\frac{1}{\lambda_m(2+\frac{\lambda_m}{M^2})}\frac{1}{(1+\frac{\lambda_m}{M^2})^{2M^2}}={\rm O}\bigl(\frac{1}{(1+\frac{1}{2})^{2M^2}}\bigr)\underset{M\to \infty}\to 0$.

As a consequence, one obtains
\[
\underset{M\to \infty}\limsup~\delta_{2}^{2}(\frac{T}{M^2})\ge \exp\bigl(-\int_{\frac{1}{2}}^{1}\frac{1}{\pi^2z^2(2+\pi^2z^2)}dz\bigr)-\exp\bigl(-\int_{\frac{1}{2}}^{1}\frac{1}{2\pi^2z^2}dz\bigr)>0.
\]

Finally,
\[
\underset{N\to\infty}\limsup~\delta_{1}^{2}(N)>0~,\quad
\underset{\Delta t\to 0}\limsup~\delta_{2}^{2}(\Delta t)>0,
\]
hence~\eqref{eq:theo_0}. This concludes the second proof of Theorem~\ref{theo:0}.

\subsection{Bounded Lipschitz test functions: proof of Theorem~\ref{theo:1}}\label{sec:proof_1}


As already explained, it is sufficient to focus on the case $r\in(\frac{1}{4},\frac{1}{2})$.

The proof is based on introducing a family $\Phi^3=\left\{\phi_{\alpha,M}^3,~\alpha\in(\frac{1}{4},\frac{1}{2}],~M\in\N\right\}\subset \mathcal{C}^{0,1}(H,\R)$, of bounded and Lipschitz continuous test functions, such that $\|\phi\|_1\le 1$ for all $\phi\in\Phi_3$. Precisely, for $\alpha\in(\frac{1}{4},\frac{1}{2}]$ and $N\in\N$, set
\[
\phi_{\alpha,M}^3(x)=\frac{\exp\bigl(-\sum_{m=M}^{\infty}\frac{|\langle x,e_m\rangle|}{\lambda_m^\alpha}\bigr)}{1+\bigl(\sum_{m=1}^{\infty}\frac{1}{\lambda_m^{2\alpha}}\bigr)^{\frac{1}{2}}}.
\]
In contrast with the families $\Phi^1$ and $\Phi^2$ introduced above, note that functions in the set $\Phi^3$ are not smooth. Introduce the notation $L_\alpha=1+\bigl(\sum_{m=1}^{\infty}\frac{1}{\lambda_m^{2\alpha}}\bigr)^{\frac{1}{2}}\in(0,\infty)$.

Let also
\[
\delta_{1}^{3}(N)=\underset{\phi\in\Phi_3}\sup |\int \phi d\mu_T-\int \phi d\mu_{T}^{(N)}|~,\quad
\delta_{2}^{3}(\Delta t)=\underset{\phi\in\Phi_3}\sup |\int \phi d\mu_T-\int \phi d\nu_{\lfloor\frac{T}{\Delta t}\rfloor}^{\Delta t}|.
\]

Let us introduce the following auxiliary function $f:[0,1]\to \R$:
\[
f(\theta)=-\log\bigl(\E[e^{-\theta |Z|}]\bigr),
\]
where $Z$ is a standard real-valued Gaussian random variable.

It is straigthforward to check that $f$ is of class $\mathcal{C}^\infty$ on $[0,1]$, that it is bounded, and that all its derivatives are bounded. Moreover, $f(0)=0$, and $f'(0)=\sqrt{\frac{2}{\pi}}$: it is crucial in the analysis below that $f'(0)\neq 0$.

\subsubsection{Spatial discretization}

First, focus on $\delta_{1}^{3}(N)$. For all $N\in\N$, and $\alpha\in(\frac{1}{4},\frac{1}{2}]$, by the indepencence property of the components of the process $X$,
\begin{gather*}
\E[\phi_{\alpha,1}^3(X(T))]=L_{\alpha}^{-1}\exp\bigl(-\sum_{n=1}^{\infty}f\bigl(\frac{\sigma_n(T)}{\sqrt{2\lambda_n}\lambda_n^\alpha}\bigr)\bigr)\\
\E[\phi_{\alpha,1}^3(X^{(N)}(T))]=L_{\alpha}^{-1}\exp\bigl(-\sum_{n=1}^{N}f\bigl(\frac{\sigma_n(T)}{\sqrt{2\lambda_n}\lambda_n^\alpha}\bigr)\bigr),
\end{gather*}
with $\sigma_n(T)^2=1-e^{-2\lambda_n T}$. Thus
\begin{align*}
\E[\phi_{\alpha,1}^{3}(X(T))]-\E[\phi_{\alpha,1}^{3}&(X^{(N)}(T))]=\E[\phi_{\alpha,1}^{3}(X(T))]\Bigl(1-\exp\bigl(\sum_{n=N+1}^{\infty}f\bigl(\frac{\sigma_n(T)}{\sqrt{2}\lambda_n^{\alpha+\frac{1}{2}}}\bigr)\bigr)\Bigr)\\
&=\E[\phi_{\alpha,1}^3(X(T))]\Bigl(1-\exp\bigl(\sum_{n=N+1}
^{\infty}\bigl(f'(0)\frac{\sigma_n(T)}{\sqrt{2}\lambda_n^{\alpha+\frac{1}{2}}}+\epsilon_n(T)\bigr)\bigr)\Bigr)
\end{align*}
where $\epsilon_n(T)=f\bigl(\frac{\sigma_n(T)}{\sqrt{2}\lambda_n^{\alpha+\frac{1}{2}}}\bigr)-f'(0)\frac{\sigma_n(T)}{\sqrt{2}\lambda_n^{\alpha+\frac{1}{2}}}={\rm O}\bigl(\frac{\sigma_n(T)^2}{2\lambda_n^{2\alpha+1}}\bigr)$.

On the one hand, $\sum_{n=N+1}^{\infty}\epsilon_n(T)\underset{N\to \infty}\to 0$. On the other hand, when $N\to \infty$, 
\begin{align*}
\sum_{n=N+1}^{\infty}\frac{\sigma_n(T)}{\sqrt{2}\lambda_n^{\alpha+\frac{1}{2}}}&=\sum_{n=N+1}^{\infty}\frac{1}{\sqrt{2}n^{2\alpha+1}\pi^{2\alpha+1}}+{\rm O}(e^{-\lambda_{N+1}T})\\
&\underset{N\to \infty}\sim \frac{C_\alpha}{N^{2\alpha}},
\end{align*}
with $C_\alpha=\int_{1}^{\infty}\frac{1}{\sqrt{2}\pi^{2\alpha+1}z^{2\alpha+1}}dz\in(0,\infty)$, by a Riemann sum argument.

Finally,
\[
\E[\phi_{\alpha,M}^3(X(T))]-\E[\phi_{\alpha,M}^3(X^{(N)}(T))]\underset{N\to \infty}\sim \frac{f'(0)\E[\phi_{\alpha,M}^3(X(T))]C_\alpha}{\lambda_N^{\alpha}}.
\]

We are now in position to conclude. Let $r\in(\frac{1}{4},\frac{1}{2})$. Then, choosing $\alpha\in(\frac{1}{4},r)$,
\[
\underset{N\to \infty}\limsup~\lambda_N^{r}\delta_{1}^{3}(N)\ge \underset{N\to \infty}\limsup~ \lambda_N^{r}\big|\E[\phi_{\alpha,1}^3(X(T))]-\E[\phi_{\alpha,1}^3(X^{(N)}(T))]\big|=\infty.
\]
This concludes the proof of Theorem~\ref{theo:1} for spatial discretization.

\subsubsection{Temporal discretization}

Now, focus on $\delta_{2}^{3}(\Delta t)$. It is convenient to choose $\Delta t=\frac{T}{M^2}$ and to consider functions $\phi_{\alpha,M}^{3}$. We claim that, for any $r\in(\frac{1}{4},\frac{1}{2})$, choosing $\alpha\in(\frac{1}{4},r)$, then 
\begin{equation}\label{eq:claim}
\underset{M\to \infty}\limsup~M^{2r}\big|\E[\phi_{\alpha,M}^3(X(T))]-\E[\phi_{\alpha,M}^{3}(X_{M^2}^{\Delta t})]\big|=\infty.
\end{equation}
On the one hand, the computations from the previous section prove that
\[
\E[\phi_{\alpha,M}^{3}(X(T))]=L_{\alpha}^{-1}\exp\bigl(-\sum_{m=M+1}
^{\infty}f\bigl(\frac{\sigma_m(T)}{\sqrt{2}\lambda_m^{\alpha+\frac{1}{2}}}\bigr)\bigr)=L_{\alpha}^{-1}\Bigl(1-\frac{f'(0)C_\alpha}{M^{2\alpha}}+{\rm O}(\frac{1}{M^{4\alpha}})\Bigr).
\]
On the other hand, using similar arguments (in particular, a Riemann sum appears),
\begin{align*}
\E[\phi_{\alpha,M}^{3}(X_{M^2}^{\Delta t})]&=L_{\alpha}^{-1}\exp\bigl(-\sum_{m=M+1}
^{\infty}f\bigl(\frac{\sigma_m(T,M)}{\lambda_m^{\alpha+\frac{1}{2}}(2+\frac{\lambda_m}{M^2})^{\frac{1}{2}}}\bigr)\bigr)\\
&=L_{\alpha}^{-1}\Bigl(1-\frac{f'(0)\overline{C}_\alpha}{M^{2\alpha}}+{\rm O}(\frac{1}{M^{4\alpha}})\Bigr),
\end{align*}
with $\sigma_m(T,M)^2=1-\frac{1}{(1+\frac{\lambda_m}{M^2})^{2M^2}}$, and
\[
\overline{C}_\alpha=\int_{1}^{\infty}\frac{1}{\sqrt{2+\pi^2z^2}\pi^{2\alpha+1}z^{2\alpha+1}}dz<C_\alpha.
\]

Thus
\[
\E[\phi_{\alpha,M}^3(X(T))]-\E[\phi_{\alpha,M}^{3}(X_{M^2}^{\Delta t})]\underset{M\to \infty}\sim \frac{f'(0)(\overline{C}_\alpha-C_\alpha)}{L_\alpha M^{2\alpha}}.
\]
This expression implies the claim~\eqref{eq:claim} holds true, hence
\[
\underset{\Delta t\to 0}\limsup~\frac{1}{\Delta t^r}\delta_{2}^{3}(\Delta t)=\infty.
\]
This concludes the proof of Theorem~\ref{theo:1} for temporal discretization.


\bibliographystyle{abbrv}
\def\cprime{$'$}


\end{document}